\input amstex
\magnification\magstephalf
\documentstyle{amsppt}

\hsize 5.72 truein
\vsize 7.9 truein
\hoffset .39 truein
\voffset .26 truein
\mathsurround 1.67pt
\parindent 20pt
\normalbaselineskip 13.8truept
\normalbaselines
\binoppenalty 10000
\relpenalty 10000
\csname nologo\endcsname 


\font\bc=cmb10
\font\tenbsy=cmbsy10

\catcode`\@=11

\def\myitem#1.{\item"(#1)."\advance\leftskip10pt\ignorespaces}

\def\qedsymbol{{\mathsurround\z@$\square$}}
\redefine\qed{\relaxnext@\ifmmode\let\next\@qed\else
  {\unskip\nobreak\hfil\penalty50\hskip2em\null\nobreak\hfil
    \qedsymbol\parfillskip\z@\finalhyphendemerits0\par}\fi\next}
\def\@qed#1$${\belowdisplayskip\z@\belowdisplayshortskip\z@
  \postdisplaypenalty\@M\relax#1
  $$\par{\lineskip\z@\baselineskip\z@\vbox to\z@{\vss\noindent\qed}}}
\outer\redefine\beginsection#1#2\par{\par\penalty-250\bigskip\vskip\parskip
  \leftline{\tenbsy x\bf#1. #2}\nobreak\smallskip\noindent}
\outer\redefine\genbeginsect#1\par{\par\penalty-250\bigskip\vskip\parskip
  \leftline{\bf#1}\nobreak\smallskip\noindent}

\def\next{\let\@sptoken= }\def\next@{ }\expandafter\next\next@
\def\@futureletnext#1{\let\nextii@#1\futurelet\next\@flti}
\def\@flti{\ifx\next\@sptoken\let\next@\@fltii\else\let\next@\nextii@\fi\next@}
\expandafter\def\expandafter\@fltii\next@{\futurelet\next\@flti}

\let\zeroindent\z@
\let\savedef@\endproclaim\let\endproclaim\relax 
\define\chkproclaim@{\add@missing\endroster\add@missing\enddefinition
  \add@missing\endproclaim
  \envir@stack\endproclaim
  \edef\endit@{\leftskip\the\leftskip\rightskip\the\rightskip}}
\let\endproclaim\savedef@
\def\thing@{.\enspace\egroup\ignorespaces}
\def\thingi@(#1){ \rm(#1)\thing@}
\def\thingii@\cite#1{ \rm\@pcite{#1}\thing@}
\def\thingiii@{\ifx\next(\let\next\thingi@
  \else\ifx\next\cite\let\next\thingii@\else\let\next\thing@\fi\fi\next}
\def\thing#1#2#3{\chkproclaim@
  \ifvmode \medbreak \else \par\nobreak\smallskip \fi
  \noindent\advance\leftskip#1
  \hskip-#1#3\bgroup\bc#2\unskip\@futureletnext\thingiii@}
\let\savedef@\endproclaim\let\endproclaim\relax 
\def\endit{\endproclaim\endit@\let\endit@\undefined}
\let\endproclaim\savedef@
\def\defn#1{\thing\parindent{Definition #1}\rm}
\def\lemma#1{\thing\parindent{Lemma #1}\sl}
\def\prop#1{\thing\parindent{Proposition #1}\sl}
\def\thm#1{\thing\parindent{Theorem #1}\sl}
\def\cor#1{\thing\parindent{Corollary #1}\sl}

\def\remk#1{\thing\zeroindent{Remark #1}\rm}
\def\example#1{\thing\zeroindent{Example #1}\rm}
\def\narrowthing#1{\chkproclaim@\medbreak\narrower\noindent
  \it\def\next{#1}\def\next@{}\ifx\next\next@\ignorespaces
  \else\bgroup\bc#1\unskip\let\next\narrowthing@\fi\next}
\def\narrowthing@{\@futureletnext\thingiii@}

\def\@cite#1,#2\end@{{\rm([\bf#1\rm],#2)}}
\def\cite#1{\in@,{#1}\ifin@\def\next{\@cite#1\end@}\else
  \relaxnext@{\rm[\bf#1\rm]}\fi\next}
\def\@pcite#1{\in@,{#1}\ifin@\def\next{\@cite#1\end@}\else
  \relaxnext@{\rm([\bf#1\rm])}\fi\next}

\advance\minaw@ 1.2\ex@
\atdef@[#1]{\ampersand@\let\@hook0\let\@twohead0\brack@i#1,\z@,}
\def\brack@{\z@}
\let\@@hook\brack@
\let\@@twohead\brack@
\def\brack@i#1,{\def\next{#1}\ifx\next\brack@
  \let\next\brack@ii
  \else \expandafter\ifx\csname @@#1\endcsname\brack@
    \expandafter\let\csname @#1\endcsname1\let\next\brack@i
    \else \Err@{Unrecognized option in @[}%
  \fi\fi\next}
\def\brack@ii{\futurelet\next\brack@iii}
\def\brack@iii{\ifx\next>\let\next\brack@gtr
  \else\ifx\next<\let\next\brack@less
    \else\relaxnext@\Err@{Only < or > may be used here}
  \fi\fi\next}
\def\brack@gtr>#1>#2>{\setboxz@h{$\m@th\ssize\;{#1}\;\;$}%
 \setbox@ne\hbox{$\m@th\ssize\;{#2}\;\;$}\setbox\tw@\hbox{$\m@th#2$}%
 \ifCD@\global\bigaw@\minCDaw@\else\global\bigaw@\minaw@\fi
 \ifdim\wdz@>\bigaw@\global\bigaw@\wdz@\fi
 \ifdim\wd@ne>\bigaw@\global\bigaw@\wd@ne\fi
 \ifCD@\enskip\fi
 \mathrel{\mathop{\hbox to\bigaw@{$\ifx\@hook1\lhook\mathrel{\mkern-9mu}\fi
  \setboxz@h{$\displaystyle-\m@th$}\ht\z@\z@
  \displaystyle\m@th\copy\z@\mkern-6mu\cleaders
  \hbox{$\displaystyle\mkern-2mu\box\z@\mkern-2mu$}\hfill
  \mkern-6mu\mathord\ifx\@twohead1\twoheadrightarrow\else\rightarrow\fi$}}%
 \ifdim\wd\tw@>\z@\limits^{#1}_{#2}\else\limits^{#1}\fi}%
 \ifCD@\enskip\fi\ampersand@}
\def\brack@less<#1<#2<{\setboxz@h{$\m@th\ssize\;\;{#1}\;$}%
 \setbox@ne\hbox{$\m@th\ssize\;\;{#2}\;$}\setbox\tw@\hbox{$\m@th#2$}%
 \ifCD@\global\bigaw@\minCDaw@\else\global\bigaw@\minaw@\fi
 \ifdim\wdz@>\bigaw@\global\bigaw@\wdz@\fi
 \ifdim\wd@ne>\bigaw@\global\bigaw@\wd@ne\fi
 \ifCD@\enskip\fi
 \mathrel{\mathop{\hbox to\bigaw@{$%
  \setboxz@h{$\displaystyle-\m@th$}\ht\z@\z@
  \displaystyle\m@th\mathord\ifx\@twohead1\twoheadleftarrow\else\leftarrow\fi
  \mkern-6mu\cleaders
  \hbox{$\displaystyle\mkern-2mu\copy\z@\mkern-2mu$}\hfill
  \mkern-6mu\box\z@\ifx\@hook1\mkern-9mu\rhook\fi$}}%
 \ifdim\wd\tw@>\z@\limits^{#1}_{#2}\else\limits^{#1}\fi}%
 \ifCD@\enskip\fi\ampersand@}

\def\pr@m@s{\ifx'\next\let\nxt\pr@@@s \else\ifx^\next\let\nxt\pr@@@t
  \else\let\nxt\egroup\fi\fi \nxt}

\define\widebar#1{\mathchoice
  {\setbox0\hbox{\mathsurround\z@$\displaystyle{#1}$}\dimen@.1\wd\z@
    \ifdim\wd\z@<.4em\relax \dimen@ -.16em\advance\dimen@.5\wd\z@ \fi
    \ifdim\wd\z@>2.5em\relax \dimen@.25em\relax \fi
    \kern\dimen@ \overline{\kern-\dimen@ \box0\kern-\dimen@}\kern\dimen@}%
  {\setbox0\hbox{\mathsurround\z@$\textstyle{#1}$}\dimen@.1\wd\z@
    \ifdim\wd\z@<.4em\relax \dimen@ -.16em\advance\dimen@.5\wd\z@ \fi
    \ifdim\wd\z@>2.5em\relax \dimen@.25em\relax \fi
    \kern\dimen@ \overline{\kern-\dimen@ \box0\kern-\dimen@}\kern\dimen@}%
  {\setbox0\hbox{\mathsurround\z@$\scriptstyle{#1}$}\dimen@.1\wd\z@
    \ifdim\wd\z@<.28em\relax \dimen@ -.112em\advance\dimen@.5\wd\z@ \fi
    \ifdim\wd\z@>1.75em\relax \dimen@.175em\relax \fi
    \kern\dimen@ \overline{\kern-\dimen@ \box0\kern-\dimen@}\kern\dimen@}%
  {\setbox0\hbox{\mathsurround\z@$\scriptscriptstyle{#1}$}\dimen@.1\wd\z@
    \ifdim\wd\z@<.2em\relax \dimen@ -.08em\advance\dimen@.5\wd\z@ \fi
    \ifdim\wd\z@>1.25em\relax \dimen@.125em\relax \fi
    \kern\dimen@ \overline{\kern-\dimen@ \box0\kern-\dimen@}\kern\dimen@}%
  }

\catcode`\@\active

\ifx\pdfoutput\undefined
  \def\pdfoutput{0}
\fi

\ifnum\pdfoutput=0
  \input epsf
\else
  \let\oldfmtname\fmtname
  \def\fmtname{plain}
  \input supp-pdf
  \let\fmtname\oldfmtname
  \define\epsfbox#1{\hbox{\convertMPtoPDF{#1}{1}{1}}}
\fi

\font\tenscr=rsfs10 
\font\sevenscr=rsfs7 
\font\fivescr=rsfs5 
\skewchar\tenscr='177 \skewchar\sevenscr='177 \skewchar\fivescr='177
\newfam\scrfam \textfont\scrfam=\tenscr \scriptfont\scrfam=\sevenscr
\scriptscriptfont\scrfam=\fivescr
\define\scr#1{{\fam\scrfam#1}}
\let\Cal\scr

\let\0\relax 
\define\restrictedto#1{\big|_{#1}}

\define\boldProj{\operatorname{\text{\bc Proj}}}
\define\boldSpec{\operatorname{\text{\bc Spec}}}
\define\Der{\operatorname{Der}}
\define\Gm{\Bbb G_{\text m}}
\define\h{{\text h}}
\define\Hom{\operatorname{Hom}}
\define\HS{\operatorname{HS}}
\define\Id{\operatorname{Id}}
\define\Image{\operatorname{Im}}
\define\Ker{\operatorname{Ker}}
\define\Spec{\operatorname{Spec}}
\define\Spf{\operatorname{Spf}}

\outer\redefine\subhead#1\par{\par\medbreak\leftline{\bc #1}\smallskip\noindent}

\topmatter
\title Jets via Hasse-Schmidt Derivations\endtitle
\author Paul Vojta\endauthor
\affil University of California, Berkeley\endaffil
\address Department of Mathematics, University of California,
  970 Evans Hall\quad\#3840, Berkeley, CA \ 94720-3840\endaddress
\date 8 January 2013 \enddate
\thanks Supported by NSF grant DMS-0200892.\endthanks

\abstract
This paper is intended to provide a general reference for jet spaces
and jet differentials, valid in maximal generality (at the level of EGA).
The approach is rather concrete, using Hasse-Schmidt (divided) higher
differentials.  Discussion of projectivized jet spaces (as in Green
and Griffiths \cite{G-G}) is included.
\endabstract
\endtopmatter

\document

This paper contains a few brief notes on how to define jets using Hasse-Schmidt
higher derivations.

I wrote them in order to better understand jets, and also to generalize the
situation more fully to the situation of arbitrary schemes.  This includes
allowing singularities, working in arbitrary characteristic
(including mixed characteristic), and working in the relative situation
of one scheme over another.

\narrowthing{}  Throughout this note, all rings (and algebras) are assumed
to be commutative.
\endit

Most of this note consists of straightforward generalizations of the
theory of derivations of order $1$, as found in Grothendieck \cite{EGA}
or Matsumura \cite{M}.  Most of this is known already to the experts,
but I am not aware of any general references, other than \cite{B-L-R,
Sect.~9.6, proof of Lemma~2}, \cite{L}, and \cite{Bu, Ch.~3, (3.14)}.
See also Green and Griffiths \cite{G-G}, in the context of complex manifolds.

Section \01 gives the basic definition of higher-order divided derivations
and differentials, leading up to the basic property of jets that they
correspond to arcs in the scheme.  The fundamental object is the algebra
$\HS^m_{B/A}$ of divided differentials,
replacing $\bigoplus_{d\ge0}S^d\Omega_{B/A}$;
here $m\in\Bbb N\cup\{\infty\}$.  The algebra $\HS^m_{B/A}$ is a graded
algebra in which the higher differentials have varying degrees; therefore
there is no obvious candidate for a module to replace $\Omega_{B/A}$.

Section \02 extends the usual first and second fundamental exact
sequences to the context of higher differentials (with some changes in the
leftmost term in each case).  Section \03 shows that
higher differentials are preserved under passing to \'etale covers; this then
implies that the definitions are preserved under localization (which can also
be shown directly).  In Section \04 this fact is used to patch together
the algebras $\HS^m_{B/A}$ to get a graded sheaf of $\Cal O_X$\snug-algebras
$\HS^m_{X/Y}$ for an arbitrary morphism $X\to Y$ of schemes.  The relative
jet space $J_m(X/Y)$ can then be defined as $\boldSpec\HS^m_{X/Y}$.
In particular, $J_m(X/Y)$ represents the functor
$Z\mapsto \Hom_Y(Z[[t]]/(t^{m+1}),X)$ from $Y$\snug-schemes $Z$ to the
category of sets (see Theorem \04.5 for notation).  If $Y=Z=\Spec k$ for
a field $k$, then this means that $k$\snug-rational points on the jet space
correspond to truncated arcs $\Spec k[t]/(t^{m+1})\to X$ over $k$.

Basic properties of this jet space are described in Section \05.
Section \06 shows that projectivized jet spaces (as in Green and
Griffiths \cite{G-G}) can be defined as $\boldProj\HS^m_{X/Y}$ for $m>0$.
Finally, Section \07 discusses how one might define log jet spaces, but any
actual work in this direction is postponed to a subsequent paper.

The Green-Griffiths projective jet spaces correspond to a quotient of
the subset of ``nonzero'' jets by the action of $\Gm$ corresponding to
linear contractions and expansions of the corresponding arcs.
Using the full automorphism group of $k[t]/(t^{m+1})$ or $k[[t]]$
gives a smaller quotient with a natural completion first defined by
J. G. Semple \cite{S} and further refined by J.-P. Demailly \cite{De}
and others.  These jet spaces may be more suitable for arithmetic applications,
since they throw out certain information which may be extraneous.
These quotients will be explored in a subsequent paper.

It should be noted that, in some fields, the term ``jet space'' has a different
meaning, namely the space parametrizing elements of the completed local
ring $\widehat{\Cal O}_{X,x}$, where $X$ is a smooth scheme over a field $k$
and $x$ is a point of $X$.  This type of jet space can be readily
defined via Grothendieck's theory of principal parts \cite{EGA 0, 20.4.14}.
See also \cite{Y}.

A more general type of (arc) jet space is described in a paper of A. Weil
\cite{W}, in which he replaces $k[t]/(t^{m+1})$ with an Artinian local ring,
in the context of real manifolds.  I thank S. Kobayashi for mentioning this
to me.

I thank Karen Smith for suggesting the use of Hasse-Schmidt derivations;
I also thank Mark Spivakovsky, Srinivas, Dan Abramovich, Mainak Poddar,
Yogesh More, Johannes Nicaise, and Martin Brandenburg for useful comments.


\beginsection{\01}{Basic definitions}

We begin with a definition.

\defn{\01.1}  Let $A$ be a ring, let $f\:A\to B$ and $A\to R$ be
$A$\snug-algebras, and let $m\in\Bbb N\cup\{\infty\}$.  Then a
{\bc higher derivation of order $m$ from $B$ to $R$ over $A$} is a
sequence $(D_0,\dots,D_m)$ (or $(D_0,D_1,\dots)$ if $m=\infty$), where
$D_0\:B\to R$ is an $A$\snug-algebra homomorphism and $D_1,\dots,D_m\:B\to R$
(or $D_1,D_2,\dots$ if $m=\infty$) are homomorphisms of (additive)
abelian groups, such that:
\roster
\myitem i.  $D_i(f(a))=0$ for all $a\in A$ and all $i=1,\dots,m$; and
\myitem ii.  (Leibniz rule) for all $x,y\in B$ and all $k=0,\dots,m$,
$$D_k(xy) = \sum_{i+j=k}D_i(x)D_j(y)\;.\tag\01.1.1$$
\endroster
\endit

\thing\zeroindent{Remarks}\rm\relax
(a).  We will often omit mention of $m$ and $A$ if they are clear
from the context.  Also, since there will never be a $D_\infty$, we will write
$(D_0,\dots,D_m)$ when $m=\infty$ to mean $(D_0,D_1,\dots)$, by abuse of
notation.

(b).  Instead of condition (i) above, Matsumura \cite{M, \S27} assumes that
the $D_i$ are all $A$\snug-module homomorphisms.  This is an equivalent
definition.  Indeed, (i) and (ii) imply that the $D_i$ are $A$\snug-module
homomorphisms.  To see the converse it suffices to show by $A$\snug-linearity
that $D_i(1)=0$ for all $i>0$; this holds by applying the Leibniz rule
to $1=1\cdot 1$.

(c).  Higher derivations were introduced by Hasse and Schmidt \cite{H-S},
and are often called {\it Hasse-Schmidt derivations.}
\endit

\example{}  If $R$ has characteristic $0$, if $R\subseteq B$, and if $D$
is a derivation from $B$ to $R$ over $A$ (in the usual sense), then
$$D_i:=\frac1{i!}D^i\;,\quad i=0,\dots,m$$
give a higher derivation of order $m$ from $B$ to $R$.
\endit

\remk{\01.2}  If $\phi\:R\to R'$ is a homomorphisms of $A$\snug-algebras
and $(D_0,\dots,D_m)$ is a higher derivation from $B$ to $R$, then
$(\phi\circ D_0,\dots,\phi\circ D_m)$ is a higher derivation from $B$ to $R'$.
Thus, if $\Der_A^m(B,R)$ denotes the set of higher derivations of order $m$
from $B$ to $R$, then $\Der_A^m(B,\cdot)$ is a covariant functor from the
category of $A$\snug-algebras to the category of sets.
\endit

\defn{\01.3}  Let $f\:A\to B$ and $m$ be as above.  Define a $B$\snug-algebra
$\HS^m_{B/A}$ to be the quotient of the polynomial algebra
$$B[x^{(i)}]_{x\in B,\;i=1,\dots,m}$$
by the ideal $I$ generated by the union of the sets
$$\gather \{(x+y)^{(i)}-x^{(i)}-y^{(i)}:x,y\in B;\; i=1,\dots,m\}\;,
    \tag\01.3.1a \\
  \{f(a)^{(i)}:a\in A;\; i=1,\dots,m\}\;,\quad\text{and} \tag\01.3.1b \\
  \Bigl\{(xy)^{(k)}-\sum_{i+j=k}x^{(i)}y^{(j)}:x,y\in B;\; k=0,\dots,m\Bigr\}\;,
    \tag\01.3.1c
  \endgather$$
where we identify $x^{(0)}$ with $x$ for all $x\in B$ (and interpret
$i=1,\dots,m$ to mean $i\in\Bbb Z_{>0}$ when $m=\infty$, by the usual
abuse of notation).  We also define the {\bc universal derivation}
$(d_0,\dots,d_m)$ from $B$ to $\HS^m_{B/A}$ by
$$d_i(x)=x^{(i)}\pmod I\;.$$
The resulting algebra $\HS^m_{B/A}$ is an algebra over $B$; it can also be
viewed as an algebra over $A$ via $f$.  It is also a graded algebra
(either over $B$ or over $A$) in which the degree of $d_ix$ is $i$.
\endit

As examples, we have
$$\align \HS^0_{B/A}&=B \\
\intertext{and}
  \HS^1_{B/A}&=\bigoplus_{d\ge0}S^d\Omega_{B/A}\tag\01.4\endalign$$
(note that $\HS^m_{B/A}$ is an algebra over $B$, not a module).

\remk{\01.5}  For $0\le i\le j\le\infty$, we have natural graded
$B$\snug-algebra homomorphisms $f_{ij}\:\HS^i_{B/A}\to\HS^j_{B/A}$.
These satisfy
$$f_{ik}=f_{jk}\circ f_{ij} \qquad\text{for all $0\le i\le j\le k\le\infty$}$$
and
$$f_{ii}=\Id\qquad\text{for all $i\in\Bbb N\cup\{\infty\}$.}$$
Thus, they form a directed system, and
$$\HS^\infty_{B/A} = \varinjlim_{m\in\Bbb N} \HS^m_{B/A}\;.\tag\01.5.1$$
\endit

\prop{\01.6}  Let $A\to B$ and $A\to R$ be $A$\snug-algebras, and let
$m\in\Bbb N\cup\{\infty\}$.  Given a derivation $(D_0,\dots,D_m)$ from
$B$ to $R$, there exists a unique $A$\snug-algebra homomorphism
$\phi\:\HS^m_{B/A}\to R$ such that
$$(D_0,\dots,D_m) = (\phi\circ d_0,\dots,\phi\circ d_m)\;.\tag\01.6.1$$
Consequently $\HS^m_{B/A}$ (together with the universal derivation) represents
the functor $\Der_A^m(B,\cdot)$.
\endit

\demo{Proof}  Define
$$\phi_0\:B[x^{(i)}]_{x\in B,\;i=1,\dots,m}\to R$$
by $x^{(i)}\mapsto D_i(x)$ for all $x\in B$, $i=0,\dots,m$, where again
$x^{(0)}$ means $x$.  The properties of a derivation imply that the kernel
of $\phi_0$ contains the ideal (\01.3.1), so we get a map $\phi$ satisfying
(\01.6.1).  Moreover, (\01.6.1) forces the choice of $\phi_0$, so $\phi$
is unique.  Thus the map
$$\Hom_A(\HS^m_{B/A},R) \to \Der_A^m(B,R)\tag\01.6.2$$
defined by the right-hand side of $(\01.6.1)$ is bijective, so the second
assertion holds.\qed
\enddemo

To relate all of the above to jets, we have the following fact.

\lemma{\01.7}  Let $f\:A\to B$, $R$, and $m$ be as above, and let
$$(D_0,\dots,D_m)\in\Der_A^m(B,R)\;.$$
Define $\phi\:B\to R[t]/(t^{m+1})$ (if $m<\infty$) or $\phi\:B\to R[[t]]$
(if $m=\infty$) by
$$x \mapsto D_0(x)+D_1(x)t+\dots+D_m(x)t^m\pmod{(t^{m+1})}\;.\tag\01.7.1$$
Then $\phi$ lies in $\Hom_A(B,R[t]/(t^{m+1}))$ or $\Hom_A(B,R[[t]])$,
respectively, and the resulting map
$$\Der_A^m(B,R) \to
  \cases \Hom_A(B,R[t]/(t^{m+1})),&\quad m<\infty; \\
    \Hom_A(B,R[[t]]),&\quad m=\infty
  \endcases\tag\01.7.2$$
is bijective.
\endit

\demo{Proof}  Since the $D_i$ are all homomorphisms of the underlying additive
groups, so is $\phi$.  The Leibniz rule implies that $\phi$ is multiplicative,
so $\phi$ is a ring homomorphism.  Finally, the conditions that $D_i(f(a))=0$
for all $a\in A$ and all $i>0$ and that $D_0$ is an $A$\snug-algebra
homomorphism imply that $\phi$ is also a homomorphism of $A$\snug-algebras.
This proves the first assertion.

The map (\01.7.2) is obviously injective, and surjectivity follows by reversing
the steps in the above paragraph to show that given a map $\phi$, (\01.7.1)
defines a higher derivation.  Thus the map (\01.7.2) is bijective.\qed
\enddemo

As an immediate consequence, we have:

\cor{\01.8} (Jet desideratum; see \cite{G, Cor.~1 of Prop.~4})
Let $A\to B$, $A\to R$, and $m$ be as in Proposition \01.6.
The bijections (\01.6.2) and (\01.7.2) define a natural bijection
$$\Hom_A(\HS^m_{B/A},R) \to
  \cases \Hom_A(B,R[t]/(t^{m+1})),&\quad m<\infty; \\
    \Hom_A(B,R[[t]]),&\quad m=\infty
  \endcases\tag\01.8.1$$
in which a map $\phi\in\Hom_A(\HS^m_{B/A},R)$ is associated to the map
$B\to R[t]/(t^{m+1})$ or $B\to R[[t]]$ given by
$$x \mapsto \phi(d_0 x) + \phi(d_1 x)t + \dots + \phi(d_m x)t^m\pmod{(t^{m+1})}
  \;.\tag\01.8.2$$
\endit

\remk{\01.9}  From now on, we will write $R[[t]]/(t^{m+1})$ instead of
$R[t]/(t^{m+1})$ (for $m<\infty$) or $R[[t]]$ (for $m=\infty$).
\endit

\remk{\01.10}  Suppose further that $R$ is an algebra over $B$; say $g\:B\to R$.
Then the bijection (\01.8.1) takes $\Hom_B(\HS^m_{B/A},R)$ to
$$\{\phi\in\Hom_A(B,R[[t]]/(t^{m+1})):z\circ\phi=g\}\;,$$
where $z\:R[[t]]/(t^{m+1})\to R$ is the map inducing the identity on $R$
and taking $t$ to $0$.
\endit

\beginsection{\02}{Fundamental Exact Sequences}

The fundamental exact sequences carry over almost directly from the case
of differentials of order $1$; the only differences stem from the fact
that we are now working with algebras instead of modules.

\thm{\02.1} (First fundamental exact sequence)  Let $A\to B\to C$ be a sequence
of ring homomorphisms, and let $m\in\Bbb N\cup\{\infty\}$.  Then the sequence
$$0 @>>> \bigl(\HS^m_{B/A}\bigr)^{+}\HS^m_{C/A} @>>> \HS^m_{C/A}
  @>>> \HS^m_{C/B} @>>> 0$$
is an exact sequence of graded $C$\snug-algebras; here $(\HS^m_{B/A})^{+}$
denotes the ideal of elements of degree $>0$ in the graded algebra
$\HS^m_{B/A}$.
\endit

\demo{Proof}  Since the term on the left is an ideal in $\HS^m_{C/A}$,
the sequence is exact on the left.

The two terms $\HS^m_{C/A}$ and $\HS^m_{C/B}$ on the right are both
quotients of the same polynomial algebra, and the generating set (\01.3.1)
for the ideal used to define $\HS^m_{C/A}$ is a subset of the corresponding
generating set for $\HS^m_{C/B}$, so the map on the right is well defined
and surjective.  The difference of the two generating sets forms a generating
set for the ideal $\bigl(\HS^m_{B/A}\bigr)^{+}$ in $\HS^m_{B/A}$, so the
sequence is exact in the middle.\qed
\enddemo

Note that the kernel of the map $\HS^m_{C/A}\to \HS^m_{C/B}$ contains
elements such as $d_1b\cdot d_1c$ with $b\in B$ and $c\in C$, so it would
not be correct to replace the term on the left with
$\bigl(\HS^m_{B/A}\bigr)^{+}\otimes_B C$.

In the case of order-1 differentials, the second fundamental exact sequence
is a refinement of the first, via the fact that $\Omega_{C/B}=0$
when $B\to C$ is surjective.  In the present case, though, $\HS^m_{C/B}=C$,
suggesting a second exact sequence of the form
$$? @>>> \bigl(\HS^m_{B/A}\bigr)^{+}\HS^m_{C/A} @>>> \HS^m_{C/A} @>>> C @>>> 0
  \;.$$
Instead, though, taking the grading into account allows us to move the $C$ two
places to the left (see also Lemma \03.4):

\thm{\02.2} (Second fundamental exact sequence)  Let $A\to B\to C$ and $m$
be as above.  Assume also that $B\to C$ is surjective, and let $I$ be
its kernel.  Let $J$ be the ideal in $\HS^m_{B/A}$,
$$J := (d_ix)_{i=0,\dots,m,\;x\in I}\;.\tag\02.2.1$$
Then the sequence
$$0 @>>> J @>>> \HS^m_{B/A} @>>> \HS^m_{C/A} @>>> 0\tag\02.2.2$$
is an exact sequence of graded $B$\snug-modules.
\endit

\demo{Proof}  Exactness on the left is obvious.  It is also obvious
from the definitions that the natural map $\alpha\:\HS^m_{B/A}\to\HS^m_{C/A}$
is surjective and that its kernel contains $J$.

From the definition of $\HS^m$, we have a commutative diagram with exact rows:
$$\CD
  0 @>>> K @>>> B[x^{(i)}]_{x\in B,\;i=1,\dots,m} @>>> \HS^m_{B/A} @>>> 0 \\
  && @VVV @VVV @VV\alpha V \\
  0 @>>> K' @>>> C[y^{(i)}]_{y\in C,\;i=1,\dots,m} @>>> \HS^m_{C/A} @>>> 0
  \endCD$$
By (\01.3.1), the leftmost vertical arrow is surjective, so by the Snake Lemma,
the kernel of the middle vertical arrow maps onto $\Ker\alpha$.  But the
kernel of the middle vertical arrow is generated by
$$I \cup \{d_ix-d_iy:i=1,\dots,m;\;x,y\in B;\; x-y\in I\}\;.$$
This implies that the kernel of $\alpha$ is generated by the set
$\{d_ix:i=0,\dots,m,\;x\in I\}$, as in (\02.2.1).\qed
\enddemo

\remk{\02.3}  In (\02.2.1), it suffices to let $x$ vary over a generating
set of $I$.
\endit

\beginsection{\03}{Formally \'Etale Algebras}

We recall the following definition:

\defn{\03.1}  Let $C$ be an algebra over $B$.  Then:
\roster
\myitem a.  $C$ is {\bc formally unramified} over $B$ if, for each
$B$\snug-algebra homomorphism $p\:D\to E$ with nilpotent kernel, and
for each $B$\snug-algebra homomorphism $v\:C\to E$, there exists
{\it at most one\/} $B$\snug-algebra homomorphism $u\:C\to D$ such that
$p\circ u=v$:
$$\epsfbox{jets1.1}$$
\myitem b.  $C$ is {\bc formally smooth} over $B$ if, for each surjective
$B$\snug-algebra homomorphism $p\:D\to E$ with nilpotent kernel, and
for each $B$\snug-algebra homomorphism $v\:C\to E$, there exists
(at least one) $B$\snug-algebra homomorphism $u\:C\to D$ such that $p\circ u=v$.
\myitem c.  $C$ is {\bc formally \'etale} over $B$ if it is formally unramified
and formally smooth over $B$.
\endroster
\endit

\remk{\03.2}  Let $(D_i)_{i\in I}$ be a directed system of $B$\snug-algebras
with maps $(f_{ij})_{i,j\in I}$, let $(p_i)_{i\in I}$ be a system of
$B$\snug-algebra homomorphisms $p_i\:D_i\to E$ compatible with the maps
$f_{ij}$, let $D=\varprojlim_{i\in I}D_i$, and let $p\:D\to E$ be the resulting
map.  If $C$ is formally unramified over $B$ and if the $p_i$ all have nilpotent
kernels, then the condition in Definition \03.1a applies also to this
map $p\:D\to E$ (i.e., there is at most one $u\:C\to D$ such that the above
diagram commutes).

Similarly, if $C$ is formally smooth over $B$, if the $p_i$ are all surjective
and have nilpotent kernels, and if the same holds for all the $f_{ij}$, then
the condition in Definition \03.1b holds for $p\:D\to E$ in this case.

In particular, these conditions apply to the map $z\:R[[t]]\to R[[t]]/(t)=R$,
so the proofs of Lemmas \03.3 and \03.5 are valid also when $m=\infty$.
\endit

\lemma{\03.3}  If $C$ is formally unramified over $B$, then $\HS^m_{C/B}=C$
for all $m\in\Bbb N\cup\{\infty\}$.
\endit

\demo{Proof}  We have the following isomorphisms:
$$ \Hom_B(\HS^m_{C/B},R) \overset\sim\to\longrightarrow
  \Hom_B(C,R[[t]]/(t^{m+1})) \overset\sim\to\longrightarrow \Hom_B(C,R)$$
for all $B$\snug-algebras $R$.  Indeed, the first arrow is Corollary \01.8.
The second arrow is defined by composing with the morphism
$z\:R[[t]]/(t^{m+1})\to R$ given by $t\mapsto 0$; it is easily seen to be
surjective; and injectivity follows by applying the definition of formal
unramifiedness to $z$, as in (\03.3.1), below.  Since this bijection holds
for all $R$, we obtain $\HS^m_{C/B}=C$.
$$\epsfbox{jets1.2}\qed\tag\03.3.1$$
\enddemo

\lemma{\03.4}  Let $A\to B\to C$ and $m$ be as in Theorem \02.1.
If $\HS^m_{C/B}=C$, then the natural map $\HS^m_{B/A}\otimes_B C\to\HS^m_{C/A}$
is surjective.
\endit

\demo{Proof}  This is a graded homomorphism of graded $C$\snug-algebras,
so it suffices to prove surjectivity on each graded part, by induction on the
degree $k$.  If $k=0$, then the map is $C\to C$, which is clearly surjective.
Let $k>0$, and assume that surjectivity has been proved in all degrees $<k$.
Then, by the first fundamental exact sequence, every homogeneous element of
degree $k$ of $\HS^m_{C/A}$ can be written as a sum of terms $\alpha\cdot\beta$,
in which $\alpha\in\HS^m_{C/A}$ and $\beta\in\bigl(\HS^m_{B/A}\bigr)^{+}$.
Moreover, we may assume that $\beta$ is of the form $d_ix$ for some
$i>0$ and $x\in B$.  But then $\alpha$ has degree $k-i<k$,
so it lies in the image of $\HS^m_{B/A}\otimes_B C$ by induction.
Thus the map is surjective in degree $k$, as was to be shown.\qed
\enddemo

\lemma{\03.5}  Let $A\to B\to C$ and $m$ be as in Theorem \02.1, and
further assume that $C$ is formally smooth over $B$.  Then the natural map
$$\alpha\:\HS^m_{B/A}\otimes_B C \to \HS^m_{C/A}$$
is invertible on the left; i.e., there exists a $C$\snug-algebra homomorphism
$$\beta\:\HS^m_{C/A}\to\HS^m_{B/A}\otimes_B C$$
such that $\beta\circ\alpha=\Id_{\HS^m_{B/A}\otimes_B C}$.
\endit

\demo{Proof}  We first claim that, for any $C$\snug-algebra $R$, the map
$$\Hom_C(\HS^m_{C/A},R) @>>> \Hom_C(\HS^m_{B/A}\otimes_B C, R)$$
is surjective.  Indeed, we have
$\Hom_C(\HS^m_{B/A}\otimes_B C,R)=\Hom_B(\HS^m_{B/A},R)$; this together
with Remark \01.10 implies that the above map is equivalent to the map
$$\split & \{\phi\in\Hom_A(C,R[[t]]/(t^{m+1})):z\circ\phi=g\} \\
  &\qquad @>>> \{\phi\in\Hom_A(B,R[[t]]/(t^{m+1})):z\circ\phi=f\}\;,
  \endsplit\tag\03.5.1$$
where $f\:B\to R$ and $g\:C\to R$ express the $B$\snug-algebra and
$C$\snug-algebra structures on $R$, respectively.  Surjectivity of (\03.5.1)
then follows from the definition of formal smoothness, as in (\03.3.1).

The lemma then follows by applying this claim to $R=\HS^m_{B/A}\otimes_B C$.\qed
\enddemo

\thm{\03.6}  Let $A\to B\to C$ be a sequence of ring homomorphisms, and
let $m\in\Bbb N\cup\{\infty\}$.  If $C$ is formally \'etale over $B$, then
the natural map
$$\HS^m_{B/A}\otimes_B C \to \HS^m_{C/A}$$
is an isomorphism of graded $C$\snug-algebras.
\endit

\demo{Proof}  It is surjective by Lemma \03.3 and injective by Lemma \03.5.\qed
\enddemo

\beginsection{\04}{Jet Schemes}

Two localization properties of the algebras $\HS^m_{B/A}$ allow one to define
$\HS^m_{X/Y}$ for an arbitrary morphism of schemes $X\to Y$.

\lemma{\04.1}  Let $f\:A\to B$ be a ring homomorphism,
and let $m\in\Bbb N\cup\{\infty\}$.  Let $S$ is a multiplicative subset of $A$.
If $f$ factors through the canonical morphism $A\to S^{-1}A$, then
the identity map on $B[x^{(i)}]_{x\in B,\;i=1,\dots,m}$ induces an
isomorphism
$$\HS^m_{B/A} \overset\sim\to\to \HS^m_{B/S^{-1}A}\;.$$
\endit

\demo{Proof}  Let $f'\:S^{-1}A\to B$ be the factored homomorphism.
The only difference between the ideals (\01.3.1) for $B/A$ and for $B/S^{-1}A$
comes from (\01.3.1b), since the image of $f'$ may be larger than the image
of $f$.  But, for all $a\in A$, $s\in S$, and $k\in\{1,\dots,m\}$, we have
$$0 = d_k f(a) = d_k(f(s)f'(s^{-1}a))
  = \sum_{i+j=k} d_i f(s)\cdot d_j f'(s^{-1}a) = f(s)d_k f'(s^{-1}a)$$
in $\HS^m_{B/A}$, by induction on $k$, so since $f(s)$ is invertible in $B$,
$d_k f'(s^{-1}a)$ vanishes in $\HS^m_{B/A}$.  Thus the two ideals
coincide.\qed
\enddemo

\lemma{\04.2}  Let $f\:A\to B$ and $m$ be as in Lemma \04.1, and let $S$
be a multiplicative subset of $B$.  Then the obvious map
$$B[x^{(i)}]_{x\in B,\;i=1,\dots,m}
  \to S^{-1}B[y^{(i)}]_{y\in S^{-1}B,\;i=1,\dots,m}$$
induces an isomorphism
$$S^{-1}\HS^m_{B/A} \overset\sim\to\to \HS^m_{S^{-1}B/A}\;.$$
\endit

\demo{Proof}  This is immediate from Theorem \03.6, since localization algebras
are formally \'etale.\qed
\enddemo

These two lemmas now allow us to define a sheaf of $\Cal O_X$\snug-algebras
for any scheme morphism $X\to Y$ (not necessarily separated).

\thm{\04.3}  Let $f\:X\to Y$ be a morphism of schemes, and let
$m\in\Bbb N\cup\{\infty\}$.  Then there exists a quasi-coherent sheaf
$\HS^m_{X/Y}$ of $\Cal O_X$\snug-algebras such that (i) for each pair of
open affines $\Spec A\subset Y$ and $\Spec B\subseteq f^{-1}(\Spec A)$,
there exists an isomorphism
$$\phi_{B/A}\:\Gamma(\Spec B,\HS^m_{X/Y})\overset\sim\to\to\HS^m_{B/A}$$
of $B$\snug-algebras, and (ii) the various $\phi_{B/A}$ are compatible
with the localization isomorphisms of Lemmas \04.1 and \04.2.
Moreover, the collection $(\HS^m_{X/Y},(\phi_{B/A})_{A,B})$ is unique up to
unique isomorphism.
\endit

\demo{Proof}  One can construct a suitable sheaf $\HS^m_{X/Y}$ by first
constructing it in the special case when $Y$ is affine, using Lemma \04.2
and the fact that the localization map of that lemma is functorial in $S$.
Then generalize to arbitrary $Y$ using both of the above lemmas, the fact
that their localization maps are functorial in the respective multiplicative
subsets, and the fact that the diagram
$$\CD T^{-1}\HS^m_{B/A} @>>> T^{-1}\HS^m_{B/S^{-1}A} \\
  @VVV @VVV \\
  \HS^m_{T^{-1}B/A} @>>> \HS^m_{T^{-1}B/S^{-1}A} \endCD$$
commutes for all $A\to B$ and all multiplicative subsets $S$ of $A$ and
$T$ of $B$ such that $A\to B$ extends to $S^{-1}A\to B$.

The existence of the $\phi_{B/A}$ follows from the above glueing process,
provided that one glues over the set of all open affines at each stage.
The uniqueness assertion follows similarly.\qed
\enddemo

\defn{\04.4}  Let $X\to Y$ and $m$ be as above.  Then the
{\bc scheme of $m$\snug-jet differentials} of $X$ over $Y$ is the scheme
$$J_m(X/Y) := \boldSpec\HS^m_{X/Y}\;.$$
If $A\to B$ is a ring homomorphism, then we also write
$$J_m(B/A)=J_m(\Spec B/\Spec A)$$
(which is equal to $\Spec\HS^m_{B/A}$).
\endit

\thm{\04.5} \cite{K-I, Prop.~2.9}
Let $X\to Y$ and $m$ be as above.  Then the scheme $J_m(X/Y)$
represents the functor from $Y$\snug-schemes to sets, given by
$$Z\mapsto \Hom_Y(Z\times_{\Bbb Z}\Spec\Bbb Z[t]/(t^{m+1}),X)\tag\04.5.1$$
if $m$ is finite, or
$$Z\mapsto \Hom_Y(Z\times_{\Bbb Z}\Spf\Bbb Z[[t]],X)\;,\tag\04.5.2$$
if $m=\infty$.
\endit

\demo{Proof}  First, for brevity, we let $Z[[t]]/(t^{m+1})$ denote the scheme
$Z\times_{\Bbb Z}\Spec\Bbb Z[t]/(t^{m+1})$ if $m$ is finite, or
the formal scheme $Z\times_{\Bbb Z}\Spf\Bbb Z[[t]]$ if $m=\infty$.
We then need to construct a bijection
$$\Hom_Y(Z[[t]]/(t^{m+1}),X) @>>> \Hom_Y(Z,J_m(X/Y))\;.\tag\04.5.3$$
Since the category of formal schemes contains the category of schemes as
a full subcategory \cite{EGA, I~10.4.8}, we will refer to $Z[[t]]/(t^{m+1})$
as a formal scheme even in the case in which $m$ is finite.  Also, we note
that if $m$ is finite and $R$ is a ring with discrete topology, then the
topology of $R[[t]]/(t^{m+1})=R[t]/(t^{m+1})$ is also discrete, so
$\Spf R[t]/(t^{m+1})=\Spec R[t]/(t^{m+1})$.  On the other hand, when $m=\infty$,
we note that the product is taken in the category of formal schemes, so
for example $(\Spec A)\times_{\Bbb Z}\Bbb Z[[t]]=\Spf A[[t]]$, by
\cite{EGA, I~10.7.2}.

The continuous ring homomorphism $\Bbb Z[[t]]/(t^{m+1})\to\Bbb Z$ given
by $t\mapsto0$ gives a morphism $Z\to Z[[t]]/(t^{m+1})$, and therefore a map
$$\Hom_Y(Z[[t]]/(t^{m+1}),X)\to\Hom_Y(Z,X)\;.$$
We also have a map $\Hom_Y(Z,J_m(X/Y))\to\Hom_Y(Z,X)$, obtained by composing
with the structural map $J_m(X/Y)\to X$.  The bijection (\04.5.3) will
be constructed so as to preserve fibers of the maps to $\Hom_Y(Z,X)$.

First consider the case in which $X$, $Y$, and $Z$ are affine, say
$X=\Spec B$, $Y=\Spec A$, and $Z=\Spec R$.  We work in the category of formal
schemes, so $A$, $B$, and $R$ will have the discrete topology and
$\Spec B=\Spf B$, etc.  Then the bijection (\04.5.3) follows immediately
from Corollary \01.8.  Indeed, by \cite{EGA, I~10.4.6 and I~10.1.3},
we have $\Hom_Y(Z[[t]]/(t^{m+1}),X)=\Hom_A(B,R[[t]]/(t^{m+1}))$.
That (\04.5.3) preserves fibers of the maps to $\Hom_Y(Z,X)$ in this case
follows by looking at the constant term in (\01.8.2).

The general case of the theorem now follows by a standard glueing argument,
working over a fixed element of $\Hom_Y(Z,X)$.\qed
\enddemo


\prop{\04.6}  For all $m$, $J_m(X/Y)$ is affine, and therefore quasi-compact
and separated, over $X$.
\endit

\demo{Proof}  Trivial.\qed
\enddemo

Let $X\to Y$ be a morphism of schemes.  Then the maps
$f_{ij}\:\HS^i_{B/A}\to\HS^j_{B/A}$ of Remark \01.5 give rise to graded
homomorphisms
$$f_{ij}\:\HS^i_{X/Y}\to\HS^j_{X/Y}\tag\04.7$$
for all $0\le i\le j\le\infty$, which again form a directed system as
in Remark \01.5.

This in turn translates into schemes: The $f_{ij}$ give rise to morphisms
$$\pi_{ji}\:J_j(X/Y)\to J_i(X/Y)\tag\04.8$$
over $X$ satisfying $\pi_{ki}=\pi_{ji}\circ\pi_{kj}$ for all
$0\le i\le j\le k\le\infty$ and $\pi_{ii}=\Id$ for all
$i\in\Bbb N\cup\{\infty\}$.

If $X\to Y$ is a morphism of schemes, then (\01.5.1) gives
$$\HS^\infty_{X/Y} = \varinjlim_{m\in\Bbb N} \HS^m_{X/Y}\tag\04.9$$
and
$$J_\infty(X/Y) = \varprojlim_{m\in\Bbb N} J_m(X/Y)\;.\tag\04.10$$

For all $m\in\Bbb N\cup\{\infty\}$, the natural graded map
$$\HS^m_{B/A}\to\HS^m_{B/B}=B$$
(as in the first fundamental exact sequence; it takes $d_jx$ to $0$ for all
$x\in B$ and all $j=1,\dots,m$) gives rise to a map
$$\HS^m_{X/Y}\to\Cal O_X\tag\04.11$$
of graded $\Cal O_X$\snug-algebras, which in turn gives rise to
``zero sections''
$$s_m\:X\to J_m(X/Y)\tag\04.12$$
satisfying $\pi_m\circ s_m=\Id_X$ for all $m$ and $\pi_{ji}\circ s_j=s_i$
for all $0\le i\le j\le\infty$.

\example{\04.13}  Let $k$ be a field of characteristic $\ne2$, let $Y=\Spec k$,
and let
$$X = \Spec k[x,y]/(y^2-x^3)$$
be the cuspidal cubic curve.  Then
$$J_1(X/Y) = \Spec k[x,y,d_1x,d_1y]/(y^2-x^3,2y\,d_1y-3x^2\,d_1x)$$
and
$$\split J_2(X/Y) = \Spec k[x,y,d_1x,d_1y,d_2x,d_2y] /& (y^2-x^3, \\
  &\;2y\,d_1y-3x^2\,d_1x, \\
  &\;(d_1y)^2+2y\,d_2y-3x(d_1x)^2-3x^2\,d_2x)\;.\endsplit$$
Therefore the fiber of $\pi_1\:J_1(X/Y)\to X$ over the point $(0,0)$
is a copy of $\Bbb A^2_k$, but the fiber of $\pi_2$ over the origin is
only isomorphic to $\Bbb A^3_{k[e]/(e^2)}$, with coordinates $d_1x$, $d_2x$,
and $d_2y$, and where $e=d_1y$.  In particular, the restriction $(d_1y)^2=0$
in the fiber over the origin occurs only in $J_2$ and not in $J_1$, so the map
$f_{12}\:\HS^1_{X/Y}\to\HS^2_{X/Y}$ is not injective, and the map
$\pi_{21}\:J_2(X/Y)\to J_1(X/Y)$ is not surjective.  See also Corollary \05.11.
\endit

\remk{\04.14}  Let $m\in\Bbb N\cup\{\infty\}$.  Then $\HS^m_{B/A}$ is
functorial in pairs $A\to B$, and $\HS^m_{X/Y}$ and $J_m(X/Y)$ are
functorial in pairs $X\to Y$.  Indeed, a commutative diagram
$$\CD B @>\phi>> B' \\
  @AAA @AAA \\
  A @>>> A' \endCD$$
of rings and ring homomorphisms induces a functorial commutative diagram
$$\CD \HS^m_{B/A} @>\HS^m_\phi>> \HS^m_{B'/A'} \\
  @AAA @AAA \\
  B @>\phi>> B'\rlap{\;,}\endCD\tag\04.14.1$$
taking $d_ib\in\HS^m_{B/A}$ to $d_i\phi(b)\in\HS^m_{B'/A'}$ for all $b\in B$
and all $i$.  This in turn induces a $B'$\snug-algebra homomorphism
$$\HS^m_{B/A}\otimes_B B'\to\HS^m_{B'/A'}\;,\tag\04.14.2$$
again functorially (via transitivity of base change for tensor).
Of particular note are the special cases $B'=B\otimes_A A'$
(base change in $A$) and $A=A'$ (functoriality in $B$).

Similarly, a commutative diagram
$$\CD X @>f>> X' \\
  @VVV @VVV \\
  Y @>>> Y'\endCD$$
of schemes induces a graded $\Cal O_{X'}$\snug-algebra homomorphism
$$\HS^m_f\:f_{*}\HS^m_{X/Y} \to \HS^m_{X'/Y'}\tag\04.14.3$$
and a commutative diagram
$$\CD J_m(X/Y) @>J_m(f)>> J_m(X'/Y') \\
  @VVV @VVV \\
  X @>f>> X'\rlap{\;,}\endCD\tag\04.14.4$$
functorially in both cases.

If $f$ is a closed immersion then so is $J_m(f)$, by the two fundamental
exact sequences (Theorems \02.1 and \02.2), or by direct computation.
Similar statements hold if $f$ is an affine morphism or a quasi-compact
morphism.

If $Y=Y'$ and $Y\to Y'$ is the identity map, then we write $J_m(f)$ as
$J_m(f/Y)$ instead.  In this case, if $f$ is an open immersion, then so
is $J_m(f/Y)$.  Consequently, if $f$ is an immersion, then so is $J_m(f/Y)$.

See also Example \05.12 for some counterexamples along these lines.
\endit

\beginsection{\05}{Functorial Properties of $\HS^m_{X/Y}$ and $J_m(X/Y)$}

\subhead Polynomial Algebras

\prop{\05.1}  Let $A$ be a ring, let $B$ be the polynomial algebra
$B=A[x_i]_{i\in I}$, and let $m\in\Bbb N\cup\{\infty\}$.  Then $\HS^m_{B/A}$
is the polynomial algebra $B[d_jx_i]_{i\in I,\;j=1,\dots,m}$.
\endit

\demo{Proof}  By Corollary \01.8,
$$\split \Hom_A(\HS^m_{B/A},R) &= \Hom_A(A[x_i]_{i\in I},R[[t]]/(t^{m+1})) \\
  &= \prod_{i\in I} R[[t]]/(t^{m+1}) \\
  &= \prod\Sb i\in I\\j=0,\dots,m\endSb R \\
  &= \Hom_A(A[x_i^{(j)}]_{i\in I,\;j=0,\dots,m},R)
  \endsplit$$
for all $A$\snug-algebras $R$.
Therefore $\HS^m_{B/A}\cong A[x_i^{(j)}]_{i\in I,\;j=0,\dots,m}$ (as an
$A$\snug-algebra).  Tracking down the above maps leads to the fact that
$x_i^{(j)}$ corresponds to $d_jx_i$ for all $i$ and $j$, so
$\HS^m_{B/A}\cong B[d_jx_i]_{i\in I,\;j=1,\dots,m}$, as a graded
$B$\snug-algebra.\qed
\enddemo

\cor{\05.2}  For any scheme $Y$,
$$\HS^m_{\Bbb A^n_Y/Y}
  \cong \Cal O_Y[x_1,\dots,x_n][d_jx_i]_{i=1,\dots,n,\;j=1,\dots,m}$$
and
$$J_m(\Bbb A^n_Y/Y) \cong \Bbb A^n_Y\times_Y\Bbb A^{nm}_Y\;.$$
\endit

\demo{Proof}  Immediate.\qed
\enddemo

\cor{\05.3}  Let $B$ be an algebra over $A$, presented as
$$B\cong A[x_i]_{i\in I}\bigm/(f_j)_{j\in J}\;.$$
Then
$$\HS^m_{B/A}
  \cong B[d_kx_i]_{i\in I,\;k=1,\dots,m}
    \bigm/(d_kf_j)_{j\in J,\;k=1,\dots,m}\;.$$
\endit

\demo{Proof}  This follows from Proposition \05.1 and the second fundamental
exact sequence (Theorem \02.2), via Remark \02.3.\qed
\enddemo

\cor{\05.4}  Let $m\in\Bbb N$.
\roster
\myitem a.  If $Y$ is a scheme and $X$ is locally of finite type over $Y$,
then $J_m(X/Y)$ is of finite type over $X$.
\myitem b.  Let
$$\CD X @>f>> X' \\
  @VVV @VVV \\
  Y @>>> Y'\endCD$$
be a commutative diagram of schemes.  If $f$ is of finite type (resp\.
locally of finite type), then so is $J_m(f)\:J_m(X/Y)\to J_m(X'/Y')$.
\endroster
\endit

\demo{Proof}  Obvious.\qed
\enddemo

\subhead Base Change

\lemma{\05.5}  Let $B$ be an $A$\snug-algebra, let $A\to A'$ be a ring
homomorphism, let $B'=B\otimes_A A'$, and let $m\in\Bbb N\cup\{\infty\}$.
Then the natural map (\04.14.2)
$$\HS^m_{B/A}\otimes_B B' \to \HS^m_{B'/A'}$$
defined by $d_ix\otimes b \mapsto bd_i(x\otimes 1)$ is an isomorphism
of graded $B'$\snug-algebras.
\endit

\demo{Proof}  As an $A'$\snug-algebra, the domain of the above map is the
same as $\HS^m_{B/A}\otimes_A A'$, and for all $A'$\snug-algebras $R$, we have
$$\split \Hom_{A'}(\HS^m_{B/A}\otimes_A A',R) &= \Hom_A(\HS^m_{B/A},R) \\
  &= \Hom_A(B,R[[t]]/(t^{m+1})) \\
  &= \Hom_{A'}(B',R[[t]]/(t^{m+1})) \\
  &= \Hom_{A'}(\HS^m_{B'/A'},R) \endsplit$$
Therefore $\HS^m_{B/A}\otimes_B B'\cong\HS^m_{B'/A'}$ as $A'$\snug-algebras.
This isomorphism is the same as the map indicated in the statement of the
lemma, so it is a $B'$\snug-algebra isomorphism, and it preserves the
grading.\qed
\enddemo

(This can also be proved explicitly by Corollary \05.3.)

Translating this into the language of schemes gives the following immediate
corollary:

\prop{\05.6}  Let $X\to Y$ be schemes, let $Y'\to Y$ be a morphism of schemes,
let $X'=X\times_Y Y'$, let $p\:X'\to X$ be the canonical projection, and
let $m\in\Bbb N\cup\{\infty\}$.  Then the homomorphism (\04.14.3)
$$p_{*}\HS^m_{X'/Y'} \to \HS^m_{X/Y}\tag\05.6.1$$
is an isomorphism of graded $\Cal O_X$\snug-algebras, and the morphism
$$J_m(X'/Y') \to J_m(X/Y)\times_X X'\tag\05.6.2$$
induced by (\04.14.4) is an isomorphism of schemes over $X'$.
\endit

\subhead Products

In this subsection, tensor products will always be taken over all $i\in I$.
Infinite tensor products are taken as in \cite{E, Prop.~16.5}.

\lemma{\05.7}  For each $i\in I$ let $B_i$ be an $A$\snug-algebra,
let $B=\bigotimes_A B_i$, and let $m$ be as above.  Then the natural map
$$\bigotimes\nolimits_B \bigl(\HS^m_{B_i/A}\otimes_{B_i} B\bigr)
  \to \HS^m_{B/A}$$
is an isomorphism of graded $B$\snug-algebras.  Here the map satisfies
$$(\dotsm\otimes 1\otimes 1\otimes (d_jx\otimes b)
    \otimes 1\otimes 1\otimes\dotsm)
  \mapsto
  b\,d_j(\dotsm\otimes 1\otimes 1\otimes x\otimes 1\otimes 1\otimes\dotsm)
  \tag\05.7.1$$
for all $b\in B$ and all $x\in B_i$, where in each expression $x$ occurs in the
$i^{\text{th}}$ place.
\endit

\demo{Proof}  As before, we first consider the corresponding $A$\snug-algebra
map
$$\bigotimes\nolimits_A \HS^m_{B_i/A} \to \HS^m_{B/A}\;.\tag\05.7.2$$
We have
$$\split \Hom_A\Bigl(\bigotimes\nolimits_A\HS^m_{B_i/A},R\Bigr)
  &= \prod_{i\in I}\Hom_A(\HS^m_{B_i/A},R) \\
  &= \prod_{i\in I}\Hom_A(B_i,R[[t]]/(t^{m+1})) \\
  &= \Hom_A(B,R[[t]]/(t^{m+1})) \\
  &= \Hom_A(\HS^m_{B/A},R) \endsplit$$
for all $A$\snug-algebras $R$, so we have an isomorphism (\05.7.2).
We leave it to the reader to check that this isomorphism satisfies (\05.7.1),
and hence is an isomorphism of graded $B$\snug-algebras.\qed
\enddemo

(Again, this result can also be proved explicitly by Corollary \05.3.)

Translating this into schemes gives as an immediate corollary (in which, again,
products are taken over all $i\in I$):

\prop{\05.8}  Let $Y$ be a scheme, let $X_i$ ($i\in I$) be $Y$\snug-schemes,
let $X=\prod_Y X_i$, and let $m$ be as above.  Then the natural map
$$\bigotimes\nolimits_{\Cal O_X}
    \bigl(\HS^m_{X_i/Y}\otimes_{\Cal O_{X_i}}\Cal O_X\bigr)
  = \bigotimes\nolimits_{\Cal O_Y} \HS^m_{X_i/Y}
  @>>> \HS^m_{X/Y}$$
is an isomorphism of graded $\Cal O_X$\snug-algebras, and
$$J_m(X/Y) \to \prod\nolimits_Y J_m(X_i/Y)
  = \prod\nolimits_X (J_m(X_i/Y)\times_{X_i} X)$$
is an isomorphism of $X$\snug-schemes.
\endit

\subhead \'Etale Morphisms

\prop{\05.9}  Let $Y$ be a scheme, let $f\:X_1\to X_2$ be an \'etale
morphism of schemes over $Y$, and let $m\in\Bbb N\cup\{\infty\}$.
Then the natural map
$$f^{*}\HS^m_{X_2/Y} \to \HS^m_{X_1/Y}$$
is an isomorphism of graded $\Cal O_{X_1}$\snug-algebras, and also
$$J_m(X_1/Y) \cong J_m(X_2/Y)\times_{X_2} X_1$$
as schemes over $X_1$.  In particular, $J_m(f/Y)$ is also \'etale.
\endit

\demo{Proof}  The first two assertions are immediate from Theorem \03.6.
The last sentence follows from the second assertion and the fact that
\'etaleness is preserved under base change.\qed
\enddemo

\subhead Smooth Morphisms

\prop{\05.10}  Let $f\:X\to Y$ be a smooth morphism and let
$m\in\Bbb N\cup\{\infty\}$.  Then $X$ is covered by open sets $U$
for which $J_m(U/Y)\cong \Bbb A^{dm}_U$, where $d\in\Bbb N$ is the
relative dimension of $U$ over $Y$.  Moreover, the sets $U$ may be chosen
independent of $m$.
\endit

\demo{Proof}  By \cite{EGA, IV~17.11.4}, smoothness of $f$ implies that $X$
is covered by open sets $U$ for which there exist integers $n\ge0$ and
\'etale morphisms $g\:U\to\Bbb P^n_Y$ such that $f\restrictedto U$
factors through $g$ and the canonical projection $\Bbb P^n_Y\to Y$.
(The converse also holds; this is a convenient criterion for smoothness.)

The result is then immediate from Corollary \05.2 and Proposition \05.9.\qed
\enddemo

\cor{\05.11}  If $f\:X\to Y$ is a smooth morphism and $0\le i\le j\le\infty$,
then $f_{ij}\:\HS^i_{X/Y}\to\HS^j_{X/Y}$ is injective and
$\pi_{ji}\:J_j(X/Y)\to J_i(X/Y)$ is surjective.  (See Example \04.13
for a counterexample when $f$ is not smooth.)
\endit

\example{\05.12}  Let $k$ be a field, let $X_1=X_2=\Bbb P^1_k$, let
$f\:X_1\to X_2$ be the map given in affine coordinates by $x\mapsto x^e$,
where $\operatorname{char} k\nmid e$, and let $m\in\Bbb Z_{>0}\cup\{\infty\}$.
Then $J_m(f/k)$ is dominant but not surjective, hence not closed.
(For $m=1$ this can be computed directly; for $m>1$, use Corollary \05.11.)
Therefore the property of being closed, proper, or projective, are not
preserved by passing from $f$ to $J_m(f/k)$ \cite{I, 2.4}.
\endit

\example{\05.13}  This example describes how a construction used in diophantine
geometry can be phrased in terms of these jet spaces.  Suppose $f\:X\to Y$
is a smooth morphism of relative dimension $n$; $\sigma\:Y\to X$ is a section;
$E$ is its image; $\Cal I$ is the ideal sheaf of $E$ (with reduced induced
subscheme structure); $U$ is some open subset of $X$ meeting $E$; and
$b\in\Cal O_X(U)$ is a function vanishing to order $\ge m$ along $E$.
After possibly shrinking $U$, we may write $\Cal I=(x_1,\dots,x_n)$ with
$x_1,\dots,x_n\in\Cal O_X(U)$.  Since $b\in\Cal I^m$, we may write
$$b = \sum\Sb \bold i\in\Bbb N^n\\|\bold i|=m\endSb
    a_{\bold i}x_1^{i_1}\dotsm x_n^{i_n}\tag\05.13.1$$
with $a_{\bold i}\in\Cal O_X(U)$ for all $\bold i$.  Now consider
$(d_mb)\restrictedto E$.  If we substitute the above expression for $b$,
and apply Leibniz' rule repeatedly, we get
$$(d_nb)\restrictedto E = \sum a_{\bold i}(d_1x_1)^{i_1}\dotsm(d_1x_n)^{i_n}$$
(all other terms involve some $x_i$ without any $d_j$ applied ($j>0$),
hence vanish when restricted to $E$).  This expression may then be regarded
as an element of $\bigl(S^m\Omega_{X/Y}\bigr)\restrictedto E$.
We need to check that it is well defined (in terms of the choice of
coefficients $a_{\bold i}$ in (\05.13.1)).  This is because, by
Corollary \05.11, the map
$$S^m\Omega_{X/Y} = \bigl(\HS^1_{X/Y}\bigr)_m @>>> \bigl(\HS^m_{X/Y}\bigr)_m$$
is injective.  Since $f$ is smooth, the above sheaves are all locally free,
so the above map remains injective when restricted to $E$.  Thus,
the construction is well defined.

Similarly, if $b$ is a local section of a line sheaf $\Cal L$ vanishing to
order $\ge m$ along $E$, then the above construction gives a local section
$(d_mb)\restrictedto E
 \in\bigl(\Cal L\otimes S^m\Omega_{X/Y}\bigr)\restrictedto E$.
\endit

\beginsection{\06}{Green-Griffiths Projective Jet Spaces}

Green and Griffiths \cite{G-G} defined projectivized jet spaces by defining
an action of $\Gm$ on $J_m(X/Y)\setminus\Image s_m$ and
constructing a quotient space $P_m(X/Y)$.  (They did the construction only
for manifolds over $\Bbb C$, but the construction readily generalizes.)
We repeat the construction here, in full generality.

\defn{\06.1}  Let $A\to B$ be a ring homomorphism and let
$m\in\Bbb N\cup\{\infty\}$.
By (\01.3.1) there is a unique graded $B$\snug-homomorphism
$$\phi_m\:\HS^m_{B/A}\to\HS^m_{B/A}[z]$$
such that $d_jx\mapsto z^jd_jx$ for all $x\in B$ and all $j=0,\dots,m$.
(This will also be denoted $\phi$ when $m$ is clear from the context.)
If $X\to Y$ is a morphism of schemes, then the $\phi_m$ over affine pieces
of $X$ and $Y$ glue together to give a unique homomorphism of graded
$\Cal O_X$\snug-algebras
$$\phi_m\:\HS^m_{X/Y}\to\HS^m_{X/Y}[z]\;,$$
and hence a unique morphism
$$\psi_m\:\Bbb A^1\times J_m(X/Y) \to J_m(X/Y)$$
of schemes over $X$.
\endit

Intuitively, one thinks of $\psi_m$ as follows.  Assume that $X$ and $Y$ are
affine, say $X=\Spec B$ and $Y=\Spec A$.  Then a closed point in $J_m(X/Y)$
corresponds to a ring homomorphism $\gamma\:B\to k[[t]]/(t^{m+1})$,
and $\psi_m(z,\gamma)$ is obtained by composing $\gamma$ with the
map $k[[t]]/(t^{m+1})\to k[[t]]/(t^{m+1})$ defined by $t\mapsto zt$.
In other words, the germ of a curve is dilated by a factor of $z$.

This can be made rigorous (via the ``jet desideratum'') as follows.
For all $A$\snug-algebras $B$ and $R$, and all $z_0\in A$, the diagram
$$\CD \Hom_A(\HS^m_{B/A},R) @>>> \Hom_A(B,R[[t]]/(t^{m+1})) \\
  @VVV @VVV \\
  \Hom_A(\HS^m_{B/A},R) @>>> \Hom_A(B,R[[t]]/(t^{m+1})) \endCD\tag\06.2$$
in which the vertical map on the left is composition with the map
$$\HS^m_{B/A} @>\phi_m>> \HS^m_{B/A}[z] @>z\mapsto z_0>> \HS^m_{B/A}$$
and the vertical map on the right is composition with
$$R[[t]]/(t^{m+1}) @>t\mapsto z_0t>> R[[t]]/(t^{m+1})\;.$$

\prop{\06.3}  Let $X\to Y$ and $m$ be as above.  Then the map $\phi_m$
satisfies:
\roster
\myitem a.  The composite map
$$J_m(X/Y) @>(0,\Id_{J_m(X/Y)})>> \Bbb A^1\times J_m(X/Y) @>\psi_m>> J_m(X/Y)$$
equals the map $s_m\circ\pi_m\:J_m(X/Y)\to J_m(X/Y)$, where
$\pi_m\:J_m(X/Y)\to X$ is the structural morphism.
\myitem b.  The composite map
$$J_m(X/Y) @>(1,\Id_{J_m(X/Y)})>> \Bbb A^1\times J_m(X/Y) @>\psi_m>> J_m(X/Y)$$
equals the identity on $J_m(X/Y)$.
\myitem c.  The maps $\psi_m$ are compatible with the maps $\pi_{ji}$:  the
diagram
$$\CD \Bbb A^1\times J_j(X/Y) @>\psi_j>> J_j(X/Y) \\
  @VV \Id_{\Bbb A^1}\times \pi_{ji} V @VV \pi_{ji} V \\
  \Bbb A^1\times J_i(X/Y) @>\psi_i>> J_i(X/Y) \endCD$$
commutes for all $0\le i\le j\le\infty$.
\myitem d.  The restriction of $\psi_m$ to
$\Gm\times(J_m(X/Y)\setminus\Image s_m)$ defines a group action, which
is faithful on the set of closed points of $J_m(X/Y)$ not lying in the
image of $s_m$.
\endroster
\endit

\demo{Proof}  This is immediate from the corresponding maps on
$\HS^m_{B/A}$.\qed
\enddemo

\defn{\06.4}  Let $X\to Y$ be a morphism of schemes and let
$m\in\Bbb Z_{>0}\cup\{\infty\}$.  Then the
{\bc Green-Griffiths projectivized jet bundle} is the scheme
$$P_m(X/Y) := \boldProj\HS^m_{X/Y}$$
over $X$.
\endit

\remk{}  If $m=0$, then the above definition yields the empty scheme.
\endit

\example{}  When $m=1$, $P_1(X/Y)\cong\Bbb P(\Omega_{X/Y})$, where
$\Bbb P(\cdot)$ is the projective space of hyperplanes, as in
\cite{EGA, II~4.1.1}.
\endit

When $m>1$, $\HS^m_{X/Y}$ is not generated over $\Cal O_X$ by its elements
of degree $1$, so $\Cal O(1)$ on $P_m(X/Y)$ is not a line sheaf.  Instead,
fibers of the structural morphism $\pi_m\:P_m(X/Y)\to X$ are
{\it weighted projective spaces\/} in the sense of Dolgacev \cite{Do}
(at least when $X$ is locally of finite type over $Y$).
For finite $m$, $\Cal O(d)$ is a line sheaf if $d$ is divisible by all
integers $1,\dots,m$.

As was the case with the earlier jet spaces, there are projection mappings
$$\pi_{ji}\:P_j(X/Y) \dashrightarrow P_i(X/Y)$$
for all $1\le i\le j\le m$, and again they form an inverse system as in (\04.8).
They are only rational maps, though (unless $i=j$).  There are also projections
$$\pi_i\:P_i(X/Y) \to X$$
(which are morphisms).  They commute with the $\pi_{ji}$:
$\pi_i\circ\pi_{ji}=\pi_j$.

\remk{\06.5}  It is easy to see that the points of
$P_m(X/Y)=\boldProj\HS^m_{X/Y}$ are the nontrivial orbits of the action
of $\Gm$, and that one can think of $P_i(X/Y)$ as consisting of certain
$\Gm$\snug-invariant prime ideals (as opposed to certain homogeneous prime
ideals).
\endit

\prop{\06.6}  These projectivized jets spaces have the following properties.
\roster
\myitem a.  (Functoriality)  A commutative diagram
$$\CD X @>f>> X' \\
  @VVV @VVV \\
  Y @>>> Y' \endCD$$
induces a commutative diagram
$$\CD U @>P_m(f)>> P_m(X'/Y') \\
  @VVV @VVV \\
  X @>f>> X'\rlap{\;,} \endCD$$
where $U$ is the open subset of $P_m(X/Y)$ corresponding to the complement
in $J_m(X/Y)$ of the pull-back via $J_m(f)$ of the image of
$s_m\:X'\to J_m(X'/Y')$.
\myitem b.  (Base change)  In the situation of part (a), if $X=X'\times_{Y'}Y$,
then $P_m(f)$ is defined everywhere on $P_m(X/Y)$ and induces an isomorphism
$$P_m(X/Y) @>\sim>> P_m(X'/Y')\times_{X'} X$$
of schemes over $X$.
\myitem c.  (Affine space)  If $X\cong Y[x_i]_{i\in I}$,
then $P_m(X/Y)\cong\Bbb P^{I,m}_X$, where $\Bbb P^{I,m}_X$ denotes
$\boldProj\Cal O_X[x_{ij}]_{i\in I,\;j=1,\dots,m}$, with $x_{ij}$ homogeneous
of degree $j$ for all $i,j$.
\myitem d.  (Closed immersions)  If a $Y$\snug-morphism $f\:X_1\to X_2$
is a closed immersion, then $P_m(f/Y)$ is defined everywhere and is a closed
immersion.
\myitem e.  (Projectivity)  If $X$ is locally of finite type over $Y$
and $m<\infty$, then $P_m(X/Y)\to X$ is projective (in the sense of
\cite{EGA, II~5.5.2}).
\myitem f.  (\'Etale morphisms)  If a $Y$\snug-morphism $f\:X_1\to X_2$
is \'etale, then $P_m(f/Y)$ is defined everywhere and induces an isomorphism
$$P_m(X_1/Y)\overset\sim\to\to P_m(X_2/Y)\times_{X_2}X_1\;.$$
\myitem g.  (Smooth morphisms)  If the morphism $X\to Y$ is smooth of relative
dimension $d$, then $X$ is covered by open sets $U$ for which
$P_m(U/Y)\cong\Bbb P^{\{1,\dots,d\},m}_U$.  Moreover, this covering
may be taken independent of $m$.
\endroster
\endit

\demo{Proof}  Left to the reader.\qed
\enddemo

\beginsection{\07}{Logarithmic jets}


As a general rule, whenever something is true in Nevanlinna theory for
holomorphic maps to complex manifolds $X$, its (appropriately stated)
counterpart for holomorphic maps to complex manifolds $X$ relative to a
normal crossings divisor $D$ also holds.  For this reason, it would be
useful to define jets in the above context.  For normal crossings divisors
on complex manifolds, the appropriate notion is jets with logarithmic poles
along the divisor (e.g., $d\log f$, $d^2\log f$, etc., where $f$ is a local
defining equation for $D$).  This concept was introduced by Noguchi \cite{N}.

I believe that this concept of log jets can also be defined in the present
context, of algebraic geometry, with arbitrary singularities allowed, but
this has not been checked, for the following reasons.

The key problem is that, when working with schemes, one uses the Zariski
topology, which is much coarser than the (classical) topology on a complex
manifold.  As an example, consider an irreducible nodal curve on a smooth
surface over a field.
$$
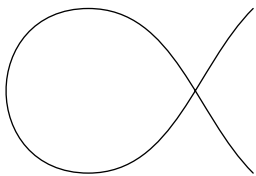
\epsfbox{jets1.3}$$
In the complex topology, it is easy to find an open neighborhood of the
node in which the two branches of the curve are defined by different
equations.  In the Zariski topology, however, this is not the case:
any open neighborhood of the intersection point will contain all but finitely
many points of the curve, and therefore the curve will always remain
irreducible.  This problem is traditionally solved by working locally in
the \'etale topology; the key result that enables one to do this seems
to be the following:\footnote{This seems to be so well known that nobody
has a reference for it.  Does anyone know of a convenient reference?}

\thm{\07.1}  Let $f\:X\to Y$ be a closed immersion of locally noetherian
excellent schemes, let $y\in Y$ be a point, let $\Cal O_y$ denote the
local ring at $y$, and let $Z$ be an irreducible component of
$X\times_Y\Spec\widehat{\Cal O}_y$.  Then there is an \'etale map $Y'\to Y$,
a point $y'\in Y'$ lying over $y$ with the same residue field $k(y')=k(y)$
(and therefore $\widehat{\Cal O}_{y'}\cong\widehat{\Cal O}_y$),
and an irreducible component $Z'$ of $X':=X\times_Y Y'$, such that
$$Z'\times_{Y'}\Spec\widehat{\Cal O}_{y'} = Z\;.\tag\07.1.1$$
\endit


\demo{Proof}  We start with some commutative algebra.  Let $A$ be an excellent
local ring with maximal ideal $\frak m$, and let $A^\h$ denote its
henselization.  By \cite{EGA, IV~18.6.6 and 18.7.6}, $A^\h$ is an excellent
local ring, and the structural map $A\to A^\h$ is a local homomorphism.
In addition, by \cite{EGA, IV~18.5.14 and 18.6.6}, the structural map
$A\to\widehat A$ factors through $A^\h$, and also $\widehat{A^\h}=\widehat A$.
If $I$ is an ideal of $A$, then by \cite{EGA, IV~18.6.8}, $(A/I)^\h=A^\h/IA^\h$.
Combining this with standard properties of completions, it follows that
the diagram
$$\CD A/I @>>> (A/I)^\h @>>> \widehat{A/I} \\
  @AAA @AAA @AAA \\
  A @>>> A^\h @>>> \widehat A\endCD$$
is Cartesian.  Finally, by \cite{EGA, proof of IV~18.6.12},
$\sqrt{IA^\h}=\sqrt IA^\h$.

Now let $A$ be the local ring $\Cal O_y$, let $I\subseteq A$ be the ideal
corresponding to $X$, and let $\sqrt I A^\h=\frak q_1\cap\dots\cap \frak q_r$
be a minimal primary decomposition of $\sqrt I A^\h=\sqrt{IA^\h}$ in $A^\h$.
It is easy to check that the $\frak q_i$ must in fact be prime.
By \cite{EGA, IV~18.9.2}, $\frak q_i\widehat A$ is also prime for all $i$.
By \cite{E, Thm.~7.2b} and \cite{M, Thm.~7.4(ii)},
$$\sqrt I\widehat A = \frak q_1\widehat A\cap\dots\cap\frak q_r\widehat A\;;$$
therefore some subset of the $\frak q_i\widehat A$ gives a primary decomposition
of $\sqrt I\widehat A$.  In particular, $Z$ comes from one of these prime
ideals, which we may assume to be $\frak q_1$.

By \cite{EGA, IV~18.6.5}, $A^\h$ is the inductive limit of all {\bc strictly
essentially \'etale} $A$\snug-algebras $B$; by definition these are
local homomorphisms $A\to B$ of local rings such that the induced map
of residue fields is an isomorphism and such that there is an \'etale
$A$\snug-algebra $C$ and an ideal $\frak n$ of $C$ lying over $\frak m$
such that $C_{\frak n}\cong B$ as $A$\snug-algebras.

Let $B$ be a strictly essentially \'etale $A$\snug-algebra containing
(finite) generating sets for each of the $\frak q_i$; then there are ideals
$\frak q_i'$ in $B$ such that $\frak q_i'A^\h=\frak q_i$ for all $i$.
By \cite{EGA, IV~18.6.5}, $B^\h=A^\h$; hence by \cite{EGA, IV~18.6.6},
$A^\h$ is faithfully flat over $B$.
Thus $\sqrt I B = \frak q_1'\cap\dots\cap\frak q_r'$.
This remains true if we replace $\frak q_i'$ with the prime ideal
$\frak q_i'A^\h\cap B$ for all $i$.

Now let $C$ be an \'etale $A$\snug-algebra and let $\frak n$ be a prime ideal
of $C$ such that $C_{\frak n}\cong B$; then there is a prime ideal $\frak p$
of $C$ corresponding to $\frak q_1$.  The generators and relations for $C$
give a scheme $Y'$ of finite type over some open affine in $Y$; let $y'$
be the point of $Y'$ corresponding to $\frak n$ and let $Z'$ be the closed
subset corresponding to $\frak p$.  Then, after replacing $Y'$ with a
suitable open subset \cite{EGA, IV~17.6.1}, it satisfies the conditions
of the theorem.\qed
\enddemo

In the slightly less general setting where $Y$ is a scheme of finite type
over a field or over an excellent Dedekind domain, this can also be handled
by applying \cite{A, Cor.~2.1}, using the radical of a product of ideals.

\defn{\07.2}  Let $X$ be a regular locally noetherian excellent scheme.
Then a {\bc normal crossings divisor} on $X$ is an effective divisor $D$ on $X$
such that, locally in the \'etale topology, irreducible components of the
support are regular and cross transversally.  If this is the case, then
we also say that $D$ {\bc has normal crossings}.
\endit

Theorem \07.1 suggests that $\HS^m_{X/Y}(\log D)$ should be defined not as
a sheaf over $X$ in the usual sense, but as an \'etale sheaf; this is
necessary so that $\HS^m_{X/Y}(\log D)$ will have the expected properties
when $Y$ is a field, $X$ is a smooth scheme over $Y$, and $D$ is a
normal crossings divisor on $X$.  Under suitably nice circumstances,
this will then descend to a quasi-coherent sheaf on the Zariski topology
(with a non-obvious definition).

In order to obtain more generality, though, it would be better to work
in the theory of logarithmic schemes.  For more details on log schemes,
see \cite{K}.

This goes far beyond the scope of the present note, so the case of log jets
will be left for future work.

\comment
||| I still need to do:

\roster
\item Jets and blowings-up:  isomorphism for jets not contained in base locus
(motivic integration lemma), $J_m$ vs. $J_{m-1}$ of blowing-up of nonsingular
subvariety.
\item approximation theorem (motivic integration seminar, 3/4).
\item Singularities and surjectivity of $\pi_{m1}$ onto the horizontal
component.
\item Crelle paper thing.
\item Closed fiber of $J_m(X/Y)$ depends only on $\Cal O_X/\frak m^{m+1}$.
\item Bijection from motivic integration (11/13).
\endroster
\endcomment


\enddocument